\pdfoutput=1
\documentclass[a4paper]{ifacconf}

\usepackage[sort&compress]{natbib}        
\usepackage{graphicx}      
\usepackage{amsmath}
\usepackage{amssymb}
\usepackage{amsfonts}
\usepackage{bbm}
\usepackage{xcolor}
\usepackage[shortlabels]{enumitem}
\usepackage{booktabs}

\newcommand{\dVolume}{{\rm d}{x}}
\newcommand{\dSurf}{{\rm d}S}
\renewcommand{\vec}[1]{\mathbf{#1}}

\newcommand{\abs}[1]{\left\vert #1 \right\vert}

\DeclareMathOperator{\trace}{trace}

\newcommand{\tr}{\mathsf{T}}
\newcommand{\oneOverG}{\frac{1}{\abs{\mathbb{G}}}}

\newcommand{\myqed}{\null\nobreak\hfill\ensuremath{\square}}

\theoremstyle{theorem}

\theoremstyle{definition}

\makeatletter
\newcommand{\subalign}[1]{%
	\vcenter{%
		\Let@ \restore@math@cr \default@tag
		\baselineskip\fontdimen10 \scriptfont\tw@
		\advance\baselineskip\fontdimen12 \scriptfont\tw@
		\lineskip\thr@@\fontdimen8 \scriptfont\thr@@
		\lineskiplimit\lineskip
		\ialign{\hfil$\m@th\scriptstyle##$&$\m@th\scriptstyle{}##$\crcr
			#1\crcr
		}%
	}
}
\newcommand*\fsize{\dimexpr\f@size pt\relax}
\makeatother

\makeatletter
\let\old@ssect\@ssect 
\makeatother
\usepackage[colorlinks=true,
			citecolor=blue,
			linkcolor=black,
			urlcolor=blue]{hyperref}
\makeatletter
\def\@ssect#1#2#3#4#5#6{%
	\NR@gettitle{#6}
	\old@ssect{#1}{#2}{#3}{#4}{#5}{#6}
}
\makeatother
\edef\endfrontmatter{%
	\unexpanded\expandafter{\endfrontmatter}
	\noexpand\endNoHyper 
}

\begin{document}
\begin{frontmatter}
\title{Verification of some functional inequalities via polynomial optimization\thanksref{footnoteinfo}}
\author{Giovanni Fantuzzi} 
\address{Department of Aeronautics, 
				Imperial College London, 
				London~SW7~2AZ 
				(e-mail: giovanni.fantuzzi10@imperial.ac.uk).\\
}
\thanks[footnoteinfo]{This work has been submitted to IFAC for possible publication. The author was funded by an Imperial College Research Fellowship.}

\begin{abstract}
	Motivated by the application of Lyapunov methods to partial differential equations (PDEs), 
	we study functional inequalities of the form $f(I_1(u),\ldots,I_k(u))\geq 0$ where $f$ is a polynomial,
	$u$ is any function satisfying prescribed constraints, 
	and $I_1(u),\ldots,I_k(u)$ are integral functionals whose integrands are polynomial in $u$, its derivatives, and the integration variable. 
	We show that such functional inequalities can be strengthened into sufficient polynomial inequalities, which in principle can be checked via semidefinite programming using standard techniques for polynomial optimization. These sufficient conditions can be used also to optimize functionals with affine dependence on tunable parameters whilst ensuring their nonnegativity.
	Our approach relies on a measure-theoretic lifting of the original functional inequality, which extends both a recent moment relaxation strategy for PDE analysis and a dual approach to inequalities for integral functionals.
\end{abstract}
\end{frontmatter}

\section{Introduction}
Lyapunov methods are a powerful tool to analyse dynamical systems governed by differential equations, wherein key system properties are inferred by constructing one or more functions of the system's state (and, possibly, of time) subject to inequality constraints. For example, classical Lyapunov functions that are positive definite and decay along trajectories prove the stability of equilibrium states~\citep{Lyapunov1892,Datko1970}. Other types of Lyapunov-like functions, generically called \textit{auxiliary functions} in this work,  can bound the effect of external disturbances~\citep{Willems1972, Dashkovskiy2013, Ahmadi2016}, certify safety~\citep{Prajna2007, Ahmadi2017, Miller2021}, approximate reachable sets and basins of attraction~\citep{Tan2006, Korda2013, Henrion2014, Valmorbida2017}, 
estimate extreme behaviour~\citep{Chernyshenko2014, Fantuzzi2016siads, Fantuzzi2020siads, Goluskin2019, Tobasco2018, Korda2021, Goluskin2020}, and solve optimal control problems~\citep{Lasserre2008, Henrion2008, Korda2018}.

For systems governed by ordinary differential equations (ODEs) with polynomial right-hand side, polynomial auxiliary functions can often be constructed via semidefinite programming. This is because the constraints on polynomial auxiliary functions reduce to polynomial inequalities imposed on the entire state space or on a subset thereof. While usually intractable, such inequalities can be strengthened by insisting that all polynomials required to be nonnegative admit a decomposition into sums of squares (SOS), which can be searched for by solving semidefinite programs (SDPs) \citep{Laurent2009,Parrilo2013,Lasserre2015}. Given enough computational resources to solve these SDPs, and barring issues with poor numerical conditioning, Lyapunov methods can therefore be applied to polynomial ODEs in a systematic way.

The use of Lyapunov methods to study dynamical systems governed by partial differential equations (PDEs), instead, remains a challenge for two reasons. One is that there is considerable freedom in the structure of auxiliary functions (more precisely, functionals) for PDEs, and it is not clear how to generalise the notion of polynomial auxiliary functions for ODEs. The second reason is that verifying if a candidate auxiliary functional satisfies all relevant constraints---such as the Lyapunov conditions for nonlinear stability---requires checking functional inequalities over infinite-dimensional function spaces. This is usually hard to do, even with computer assistance.

\cite{Valmorbida2015c,Valmorbida2015b} and \cite{Ahmadi2016,Ahmadi2017,Ahmadi2019} demonstrated that SOS polynomials and SDPs can be used to verify the nonnegativity of integral functionals with polynomial integrands in one or two spatial dimensions. This is achieved without discretization of the PDE state, but rather by requiring the polynomial integrand to be nonnegative pointwise after augmenting it by terms that integrate to zero. This approach, generalized to arbitrary spatial dimensions by~\cite{Chernyavsky2021}, enables one to search for auxiliary functions in the form of \textit{integrals of polynomials} of the PDE state. A dual formulation of the same strategy, based on so-called \textit{occupation measures}, was independently introduced by \cite{Korda2018} to relax optimal PDE control problems into SDPs.

Unfortunately, auxiliary functionals for PDEs that can be represented as integrals of polynomials may not be sufficiently general. For example, Lyapunov analysis of fluid flows \citep{Goulart2012,Huang2015,Goluskin2019,Fuentes2019} required functionals of the form
\begin{equation}\label{e:general-aux-function}
\mathcal{V}(u) = V(\mathcal{I}_1(u),\ldots,\mathcal{I}_k(u)),
\end{equation}
where  $V$ is a $k$-variate polynomial and $\mathcal{I}_1(u),\ldots,\mathcal{I}_k(u)$ are fixed integral functionals of the fluid's velocity  $u$. Such auxiliary functionals were optimized by projecting the governing PDEs onto a finite-dimensional ODE and estimating the projection error.

Motivated by the desire to avoid application-specific estimates, here we study the fundamental problem of verifying if functionals in the form~\eqref{e:general-aux-function} are nonnegative for all $u$ satisfying prescribed constraints such as boundary conditions. Leveraging the occupation measure framework of \cite{Korda2018}, we show that the broad class of functional inequalities obtained when~\eqref{e:general-aux-function} is a \emph{polynomial of integrals} can be strengthened into sufficient polynomial inequalities over finite-dimensional sets specified by linear constraints and linear matrix inequalities. These purely algebraic problems can be further strengthened into SDPs using standard techniques for polynomial optimization, so one can attempt to verify hard functional inequalities using algorithms for convex programming. In principle, therefore, one can implement Lyapunov methods for PDEs with no need for rigorous ODE projections. Doing this in practice, of course, requires overcoming well-known challenges related to the numerical solution of SDPs, which may be too large or too ill-conditioned for currently available general-purpose solvers. We do not address such implementation issues here and leave them to future work. Nevertheless, we demonstrate that our approach to verifying functional inequalities yields sharp results on two relatively simple examples. 

\section{A class of functional inequalities}\label{s:problem-statement}
Let $\Omega \subset \mathbb{R}^n$ be an open bounded domain with smooth boundary $\partial\Omega$. 
Given an integer $p$ with $1<p<\infty$, let $W^{1,p}(\Omega;\mathbb{R}^m)$ be  the Sobolev space of $p$-integrable functions $u:\Omega \to \mathbb{R}^m$ with $p$-integrable gradients. The $m \times n$ Jacobian matrix $\nabla u$ has entries $(\nabla u)_{ij} = \frac{\partial u_i}{\partial x_j}$.

Given an integer $d$, let $\mathcal{S}$ be a subset of the multi-index set
\begin{align}
\!\mathcal{E}_{d,p} := \{ 
(\alpha, \beta,\gamma) \in \mathbb{N}^n \times \mathbb{N}^m \times \mathbb{N}^{m \times n}\!:
&\abs{\alpha} \leq d,\nonumber\\
&\abs{\beta} + \abs{\gamma} \leq p \}.
\label{e:exp-set}
\end{align} 
Consider a functional $\mathcal{F}$ on $W^{1,p}(\Omega;\mathbb{R}^m)$ defined as
\begin{equation}
\mathcal{F}(u) := 
f \bigg( \left\{ \int_\Omega x^\alpha u(x)^\beta \nabla u(x)^\gamma \,\,\dVolume \right\}_{(\alpha,\beta,\gamma) \in \mathcal{S}} \bigg),
\label{e:polynomial-of-integrals}
\end{equation}
where $f$ is a function in $\abs{\mathcal{S}}$ variables. Here and throughout the paper, $x^\alpha = x_1^{\alpha_1} \cdots x_n^{\alpha_n}$ is the multivariate monomial with variables $x = (x_1,\ldots,x_n)$, exponent $\alpha \in \mathbb{N}^n$, and degree $\abs{\alpha} = \alpha_1 + \cdots +\alpha_n$. The multivariate monomials $u^\beta$ and $\nabla u^\gamma$ are defined similarly. The bound $\abs{\beta} + \abs{\gamma} \leq p$ in~\eqref{e:exp-set} ensures that the integrals in~\eqref{e:polynomial-of-integrals} are well-defined. 

We are interested in the following fundamental problems.

\begin{prob}[Feasibility]\label{p:main-problem}
	Check if $\mathcal{F}(u)$ in~\eqref{e:polynomial-of-integrals} is nonnegative for all functions $u$ in the subset $\mathcal{U}$ of $W^{1,p}(\Omega;\mathbb{R}^m)$ defined by the differential and boundary conditions
	%
	\begin{subequations}
		\label{e:pde-bcs}
		\begin{align}
		\label{e:pde-bcs-a}
		a\left({x},{u}({x}),\nabla{u}({x})\right) &= 0 \quad \text{on } \Omega,\\
		b\left({x},{u}({x})\right) &= 0 \quad \text{on } \partial\Omega,
		\label{e:pde-bcs-b}
		\end{align}
	\end{subequations}
	where
	$a:\overline{\Omega} \times \mathbb{R}^m \times \mathbb{R}^{m \times n} \to \mathbb{R}$ and
	$b:\partial\Omega \times \mathbb{R}^m \to \mathbb{R}$
	are given continuous functions that do not grow faster than a degree-$p$ polynomial in the second and third arguments.
\end{prob}

\begin{prob}[Optimization]\label{p:main-opt-problem}
	Suppose $f$ in~\eqref{e:polynomial-of-integrals} depends affinely on a vector of parameters $\lambda \in \mathbb{R}^\ell$ and let $\mathcal{U}$ be as in Problem~\ref{p:main-problem}. Given a convex function $b:\mathbb{R}^\ell \to \mathbb{R}$, solve the convex optimization problem
	\begin{equation}\label{e:opt-problem-general}
	\min_{\lambda \in \mathbb{R}^\ell} \, b(\lambda): \; \mathcal{F}(u) \geq 0 \quad \forall u \in \mathcal{U}.
	\end{equation}
\end{prob}

Our approach to these problems is to derive sufficient but more tractable conditions for the functional inequality $\mathcal{F}(u)\geq 0$. We do so in three steps. First, we lift the functional inequality into an equivalent inequality on so-called \emph{occupation} and \emph{boundary measures} \citep{Korda2018}. We then strengthen the latter by enforcing it over a larger convex set of measures, obtaining a measure-theoretic condition whose validity suffices to conclude that $\mathcal{F}(u)\geq 0$ on $\mathcal{U}$. Finally, for problems with polynomial data (cf.~\S\ref{s:sdp-projection}) we replace this measure-theoretic condition with a sufficient constrained polynomial inequality that, at least in principle, can be tested computationally. This polynomial inequality depends affinely on $\lambda$ when so does $\mathcal{F}$, so it can be used to strengthen the constraint in~\eqref{e:opt-problem-general} and compute a feasible (and, sometimes, near-optimal) $\lambda$. 


\begin{rem}
	Our results immediately extend to functions $f$ in~\eqref{e:polynomial-of-integrals} that depend also on boundary integrals. 
	Moreover, it is straightforward to consider function sets $\mathcal{U}$ whose definitions, in addition to the equalities in~\eqref{e:pde-bcs}, include differential and boundary inequalities, equality or inequality constraints on integrals over the domain or the boundary, and multiple constraints of each type. Finally, only minor modifications are needed to 
	consider domains $\Omega$ whose boundary is only piecewise-smooth and Lipschitz.
\end{rem}

\begin{rem}
	If the function $f$ in~\eqref{e:polynomial-of-integrals} is linear, then $\mathcal{F}$ is an integral functional and may be rewritten as
	\begin{equation*}
	\mathcal{F}(u) = \int_\Omega \phi({x},{u}(x),\nabla {u}(x)) \,\,\dVolume
	\end{equation*}
	for some polynomial $\phi$. In such cases, a duality argument \citep{Fantuzzi2019arxiv} shows that our sufficient conditions for the functional inequality in Problem~\ref{p:main-problem} are equivalent to those proposed by~\cite{Valmorbida2015c,Valmorbida2015b}, \cite{Ahmadi2016,Ahmadi2017,Ahmadi2019}, and \cite{Chernyavsky2021}. 
\end{rem}

\section{Lifting via occupation measures}
\label{s:om-lifting}

This section describes the first two steps of the strategy outlined above. Occupation and boundary measures are briefly reviewed in~\S\ref{s:ob-measures}, while sufficient measure-theoretic conditions for the inequality $\mathcal{F}(u)\geq 0$ on $\mathcal{U}$ are derived in~\S\ref{s:suff-condition-measures}. 
To lighten the notation, we will often write $\langle f, \mu \rangle$ to indicate the integral of a function $f$ from a set $\Theta$ into $\mathbb{R}$ against a measure $\mu$ supported on $\Theta$, i.e., 
\begin{equation*}
\langle f, \mu \rangle := \int_{\Theta} f(\vartheta) \, {\rm d}\mu.
\end{equation*}

\subsection{Occupation and boundary measures}\label{s:ob-measures}
Given any function $u \in W^{1,p}(\Omega;\mathbb{R}^m)$, let $\mu_u$ be the pushforward of the Lebesgue measure on $\Omega$ through the map $x \mapsto (x, u(x), \nabla u(x))$. By construction, $\mu_u$ is supported on the product set
\begin{equation}\label{e:D-space}
D := \overline{\Omega} \times \mathbb{R}^{m} \times \mathbb{R}^{m\times n}
\end{equation}
and satisfies
\begin{equation}\label{e:occ-measure-def}
\int_D \eta(x,y,Z) \,{\rm d}\mu_u
= \int_{\Omega} \eta(x,u(x),\nabla u(x)) \,\,\dVolume
\end{equation}
for all continuous functions $\eta:D \to \mathbb{R}$ that grow no faster than a degree-$p$ polynomial in the last two arguments (so the right-hand side is well-defined). Following \cite{Korda2018}, we call $\mu_u$ the \emph{occupation measure} generated by $u$.
In complete analogy, the \emph{boundary measure} $\nu_u$ generated by $u \in W^{1,p}(\Omega;\mathbb{R}^m)$ is the pushforward of the surface measure on $\partial\Omega$ through the map $x \mapsto (x, u(x))$. By construction, $\nu_u$ is supported on the product set
\begin{equation}\label{e:B-space}
B : =\partial \Omega \times \mathbb{R}^{m}
\end{equation}
and satisfies
\begin{equation}\label{e:bnd-measure-def}
\int_B \zeta(x,y) \,{\rm d}\nu_u
= \int_{\partial\Omega} \zeta(x,u(x)) \,\dSurf
\end{equation}
for all continuous functions $\zeta:B \to \mathbb{R}$ that grow no faster than a degree-$p$ polynomial in the second argument.

The occupation and boundary measures $\mu_u$ and $\nu_u$ are related via identities implied by the divergence theorem. Let $\Phi$ be the linear space of continuously differentiable $n$-valued functions $\varphi: \overline{\Omega}\times \mathbb{R}^m \to \mathbb{R}^n$ such that, given any $u \in W^{1,p}(\Omega;\mathbb{R}^m)$, the function $x \mapsto \varphi(x,u(x))$ is differentiable using the chain rule and has an integrable derivative. To every $\varphi \in \Phi$ we associate the ``total divergence'' function
\begin{equation*}
\mathcal{D}\varphi(x,y,Z) = \sum_{i=1}^n \bigg[ 
	\frac{\partial\varphi_i}{\partial x_i} (x,y) + 
	\sum_{j=1}^m Z_{ji} \frac{\partial\varphi_i}{\partial y_j} (x,y) \bigg]
\end{equation*}
from $\overline{\Omega}\times \mathbb{R}^m \times \mathbb{R}^{m \times n}$ into $\mathbb{R}$. This definition ensures that $\mathcal{D}\varphi(x,u(x),\nabla u(x))$ is the divergence of $\varphi(x,u(x))$ when the latter is calculated using the chain rule. With this notation, it is straightforward to combine the divergence theorem with identities~\eqref{e:occ-measure-def} and~\eqref{e:bnd-measure-def} to conclude the following result. (See also Theorem~1 in \citealp{Korda2018}.)

\begin{lem}\label{lemma:div-theorem}
	Let $\mu_u$ and $\nu_u$ be the occupation and boundary measures generated by a function $u \in W^{1,p}(\Omega;\mathbb{R}^m)$, and let $\hat{n}(x)$ be the outwards unit normal vector to $\partial\Omega$. Then,
	\begin{equation}\label{e:measure-ibp}
	\langle \mathcal{D}\varphi, \mu_u  \rangle = \langle \varphi \cdot \hat{n} , \nu_u \rangle \quad \forall \varphi \in \Phi.
	\end{equation}
\end{lem}

\subsection{Measure-theoretic strengthening of Problem~\ref{p:main-problem}}\label{s:suff-condition-measures}

Consider Problem~\ref{p:main-problem}, which asks to verify if $\mathcal{F}(u) \geq 0$ for all functions $u$ in the set $\mathcal{U}$ defined by the differential and boundary constraints in~\eqref{e:pde-bcs}. Since identity~\eqref{e:occ-measure-def} implies that
\begin{equation*}
\mathcal{F}(u) = f(\{\langle x^\alpha y^\beta Z^\gamma, \mu_u \rangle\}_{(\alpha,\beta,\gamma) \in \mathcal{S}}),
\end{equation*}
it is natural to extend $\mathcal{F}$ into a functional on a suitable space of measures, and then reformulate the functional inequality $\mathcal{F}(u)\geq 0$ as a constraint on occupation measures. 

To do this, let $\mathcal{M}_+(D)$ and $\mathcal{M}_+(B)$ be the cones of Borel measures supported on the sets $D$ and $B$ defined in~\eqref{e:D-space} and~\eqref{e:B-space}, respectively. The functional
\begin{equation}\label{e:G-def}
\mathcal{G}(\mu,\nu) = f(\{\langle x^\alpha y^\beta Z^\gamma, \mu \rangle\}_{(\alpha,\beta,\gamma) \in \mathcal{S}})
\end{equation}
from $\mathcal{M}_+(D) \times \mathcal{M}_+(B)$ into $\mathbb{R}$ clearly satisfies $\mathcal{G}(\mu_u,\nu_u) = \mathcal{F}(u)$ whenever $(\mu_u,\nu_u)$ is a pair of occupation and boundary measures generated by a function $u \in W^{1,p}(\Omega;\mathbb{R}^m)$. Problem~\ref{p:main-problem} may therefore be restated as follows.

\begin{prob}\label{p:measure-theoretic-problem}
	Check if the functional $\mathcal{G}(\mu,\nu)$ in~\eqref{e:G-def} is nonnegative when $(\mu,\nu) \in \mathcal{M}_+(D) \times \mathcal{M}_+(B)$ is a pair of occupation and boundary measures generated by $u \in \mathcal{U}$.
\end{prob}

\begin{rem}
	The apparently redundant dependence of $\mathcal{G}$ on $\nu$ is useful in what follows to exploit the relations between occupation and boundary measures given by Lemma~\ref{lemma:div-theorem}. Moreover, it enables us to easily extend our approach to functionals $\mathcal{F}(u)$ that depend also on the surface integrals $\int_{\partial\Omega} x^\alpha u(x)^\beta \,\dSurf$ with $\abs{\alpha}\leq d$ and $\abs{\beta}\leq p$, which can be rewritten as $\langle x^\alpha y^\beta, \nu_u \rangle$ using~\eqref{e:bnd-measure-def}.
\end{rem}

Problem~\ref{p:measure-theoretic-problem} is just as hard as verifying that $\mathcal{F}(u) \geq 0$ on the function set $\mathcal{U}$. However, our measure-theoretic reformulation suggests that to check this functional inequality it is enough to verify that $\mathcal{G}$ is nonnegative on a subset $\mathcal{A} \subseteq \mathcal{M}_+(D) \times \mathcal{M}_+(B)$ that includes all occupation and boundary measures generated by $u \in \mathcal{U}$. For suitably chosen $\mathcal{A}$, this sufficient measure-theoretic condition can be more tractable that the original functional inequality.

One particularly convenient strategy is to define $\mathcal{A}$ using linear constraints necessarily satisfied by occupation and boundary measures generated by $u \in \mathcal{U}$. Lemma~\ref{lemma:div-theorem} already provides an uncountable number of such constraints. Additional ones can be derived from identities~\eqref{e:occ-measure-def} and~\eqref{e:bnd-measure-def}, and from the differential and boundary conditions~\eqref{e:pde-bcs-a} and~\eqref{e:pde-bcs-b} imposed on $u \in \mathcal{U}$. Specifically, let $a$ and $b$ be the functions appearing in~\eqref{e:pde-bcs}. Further, let $\Psi_a$ and $\Psi_b$ be the sets of continuous functions $\psi:D\to \mathbb{R}$ and $\rho:B \to \mathbb{R}$ such that the functions $(\psi a)(x,y,Z)$ and $(\rho b)(x,y)$ grow no faster than a degree-$p$ polynomial in $y$ and $Z$. Then,

\begin{lem}\label{lemma:A-def}
	Every pair of occupation and boundary measures generated by $u \in \mathcal{U}$ belongs to the set $\mathcal{A}$ of pairs $(\mu,\nu) \in \mathcal{M}_+(D) \times \mathcal{M}_+(B)$ satisfying
	\begin{subequations}
		\begin{align}
		\label{e:meas-cond-div}
		\langle \mathcal{D}\varphi, \mu  \rangle &= \langle \varphi \cdot \hat{n} , \nu \rangle &&\forall \varphi \in \Phi,\\
		\label{e:meas-cond-marginal-mu}
		\langle \eta, \mu \rangle &= \textstyle\int_{\Omega} \eta(x) \,\,\dVolume  &&\forall \eta \in C(\overline{\Omega}),\\
		\label{e:meas-cond-marginal-nu}
		\langle \zeta, \nu \rangle &= \textstyle\int_{\partial\Omega} \zeta(x) \,\dSurf &&\forall \zeta \in C(\partial\Omega),\\
		\label{e:meas-cond-diff}
		\langle \psi a, \mu \rangle &= 0 &&\forall \psi \in \Psi_a,\\
		\label{e:meas-cond-bc}
		\langle \rho b, \nu \rangle &= 0 &&\forall \rho \in \Psi_b.
		\end{align}
	\end{subequations}
\end{lem} 
\begin{pf}
	If $\mu$ and $\nu$ are occupation and boundary measures generated by $u \in \mathcal{U}$, condition~\eqref{e:meas-cond-div} is just a restatement of Lemma~\ref{lemma:div-theorem}. Conditions~\eqref{e:meas-cond-marginal-mu} and~\eqref{e:meas-cond-marginal-nu} follow directly from identities~\eqref{e:occ-measure-def} and~\eqref{e:bnd-measure-def}. Condition~\eqref{e:meas-cond-diff} holds because
	\begin{multline*}
	\langle \psi a, \mu \rangle 
	= \int_D \psi(x,y,Z)a(x,y,Z) \,{\rm d}\mu \nonumber\\
	= \int_\Omega \psi(x,u,\nabla u) \, a(x,u,\nabla u) \,\,\dVolume
	= 0,
	\end{multline*}
	where the last two equalities are a consequence of~\eqref{e:occ-measure-def} and~\eqref{e:pde-bcs-a}, respectively. Condition~\eqref{e:meas-cond-bc} is proven analogously using~\eqref{e:bnd-measure-def} and~\eqref{e:pde-bcs-b}.
	\myqed
\end{pf}

Combining Lemma~\ref{lemma:A-def} with the preceding discussion yields the following result. 

\begin{prop}\label{prop:sufficient-measure-condition}
	Let  $\mathcal{G}(\mu,\nu)$ be as in~\eqref{e:G-def}, and let $\mathcal{A}$ be the set of measures $(\mu,\nu)\in\mathcal{M}_+(D) \times \mathcal{M}_+(B)$ satisfying~\eqref{e:meas-cond-div}--\eqref{e:meas-cond-bc}. If
	\begin{equation}\label{e:G-suff-condition}
	\mathcal{G}(\mu,\nu) \geq 0 \quad \forall (\mu,\nu) \in \mathcal{A},
	\end{equation}
	then the functional $\mathcal{F}(u)$ in~\eqref{e:polynomial-of-integrals} is nonnegative for all functions $u$ in the set $\mathcal{U}$ defined in Problem~\ref{p:main-problem}.
\end{prop}

\begin{rem}\label{remark:dirichlet-bcs}
	If the functions $u \in \mathcal{U}$ satisfy Dirichlet boundary conditions, meaning that $b(x,u) = u - h(x)$ in~\eqref{e:pde-bcs-b} for some fixed function $h$, then their boundary measures $\nu_u$ coincide with the known measure $\nu_h$. In such cases, one can replace $\nu$ with $\nu_h$ in Lemma~\ref{lemma:A-def} and Proposition~\ref{prop:sufficient-measure-condition}. Then, condition~\eqref{e:meas-cond-marginal-nu} can be dropped, while the term $\langle \varphi \cdot \hat{n} , \nu \rangle$ in~\eqref{e:meas-cond-div} becomes the real number $\langle \varphi \cdot \hat{n} , \nu_h \rangle = \int_{\partial\Omega}\varphi(x,h(x)) \cdot \hat{n}(x) \,\dSurf$, which in principle can be computed explicitly for every choice of $\varphi \in \Phi$.
\end{rem}

\subsection{Exploiting symmetries}\label{s:symmetries}

Proposition~\ref{prop:sufficient-measure-condition} can be refined if the functional inequality in Problem~\ref{p:main-problem} enjoys symmetries. Specifically, if the domain $\Omega$, the integrals 
entering~\eqref{e:polynomial-of-integrals}, and the constraints in~\eqref{e:pde-bcs} are all invariant under a group of linear transformations, then it suffices to verify~\eqref{e:G-suff-condition} on the group-invariant subset of~$\mathcal{A}$. 

To make these idea precise, let $\mathbb{G}$ be a group of orthonormal linear transformations on $\mathbb{R}^n \times \mathbb{R}^m$, meaning a family (finite, countable or uncountable) of pairs of matrices $(A,B) \in \mathbb{R}^{n \times n} \times \mathbb{R}^{m \times m}$ such that:
\begin{enumerate}[({G}1),  leftmargin=*]
	\item\label{ass:group-1} $A^\tr = A^{-1}$ and $B^\tr = B^{-1}$;
	\item\label{ass:group-2} If $(A,B) \in \mathbb{G}$, then $(A^{-1},B^{-1}) \in \mathbb{G}$;
	\item\label{ass:group-3} If $(A_1,B_1),\, (A_2,B_2) \in \mathbb{G}$, then $(A_1A_2,B_1B_2) \in \mathbb{G}$.
\end{enumerate}
Given $(A,B) \in \mathbb{G}$, define linear operators $G_{A,B}$ and $H_{A,B}$ on $\mathbb{R}^n \times \mathbb{R}^m \times \mathbb{R}^{m\times n}$ and $\mathbb{R}^n \times \mathbb{R}^m$, respectively, via
\begin{subequations}
	\begin{align*}
	G_{A,B}(x,y,Z) &= \left( Ax, \, By,\, B Z A^\tr \right)\\
	H_{A,B}(x,y) &= \left( Ax, \, By \right).
	\end{align*}
\end{subequations}
%
We say that Problem~\ref{p:main-problem} is invariant under $\mathbb{G}$ if:
\begin{enumerate}[({A}1), leftmargin=*]
	\item\label{ass:invariance-1} $x \in \Omega \iff Ax \in \Omega$ for all $(A,B) \in \mathbb{G}$;
	\item\label{ass:invariance-2} $a \circ G_{A,B} = a$ for all $(A,B) \in \mathbb{G}$;
	\item\label{ass:invariance-3} $b \circ H_{A,B} = b$ for all $(A,B) \in \mathbb{G}$;
	\item\label{ass:invariance-4} $(Ax)^\alpha (By)^\beta (B Z A^\tr)^\gamma = x^\alpha y^\beta Z^\gamma$ for all $(\alpha,\beta,\gamma) \in \mathcal{S}$ and all $(A,B) \in \mathbb{G}$.
\end{enumerate}


In such cases, Proposition~\ref{prop:sufficient-measure-condition} can be replaced with
\begin{prop}\label{prop:sufficient-measure-condition-symm}
	Let  $\mathcal{G}(\mu,\nu)$ be as in~\eqref{e:G-def}, and let $\mathcal{A}_{\mathbb{G}}$ be the set of measures $(\mu,\nu) \in \mathcal{M}_+(D) \times \mathcal{M}_+(B)$ that satisfy~\eqref{e:meas-cond-div}--\eqref{e:meas-cond-bc} as well as the invariance conditions
	\begin{subequations}
		\label{e:invariance}
		\begin{align}
		\label{e:invariance-mu}
		\langle \eta \circ G_{A,B}, \mu \rangle &= \langle \eta, \mu \rangle &&\forall\eta \in C(D),\, (A,B) \in \mathbb{G},\\
		\label{e:invariance-nu}
		\langle \zeta \circ H_{A,B}, \nu \rangle &= \langle \zeta, \nu \rangle &&\forall\zeta \in C(B),\, (A,B) \in \mathbb{G}.
		\end{align}
	\end{subequations}
	If assumptions~\ref{ass:invariance-1}--\ref{ass:invariance-4} hold and 
	\begin{equation}\label{e:G-suff-condition-symm}
	\mathcal{G}(\mu,\nu) \geq 0 \quad \forall (\mu,\nu) \in \mathcal{A}_{\mathbb{G}},
	\end{equation}
	then the functional $\mathcal{F}(u)$ in~\eqref{e:polynomial-of-integrals} is nonnegative for all functions $u$ in the set $\mathcal{U}$ defined in Problem~\ref{p:main-problem}.
\end{prop}
	
This result is a direct consequence of the following lemma. 
\begin{lem}
Under conditions~\ref{ass:invariance-1}--\ref{ass:invariance-4}, for every $(\mu,\nu) \in \mathcal{A}$ there exists $(\overline{\mu},\overline{\nu}) \in \mathcal{A}_{\mathbb{G}}$ such that $\mathcal{G}(\overline{\mu},\overline{\nu}) = \mathcal{G}(\mu,\nu)$.
\end{lem}
\begin{pf}
	Let $\abs{\mathbb{G}}$ be the cardinality of $\mathbb{G}$. Fix $(\mu,\nu) \in \mathcal{A}$ and define the group-averaged measures
	\begin{subequations}
	\begin{align*}
	\overline{\mu}&:= \oneOverG\sum_{(A,B) \in \mathbb{G}} \!\!G_{A,B} \sharp \mu, \\[-0.35ex]
	\overline{\nu}&:= \oneOverG\sum_{(A,B) \in \mathbb{G}} \!\!H_{A,B} \sharp \nu,
	\end{align*}
	\end{subequations}
	where $G_{A,B} \sharp \mu$ and $H_{A,B} \sharp \nu$ are the pushforwards of $\mu$ and $\nu$ by $G_{A,B}$ and $H_{A,B}$. The identity $\mathcal{G}(\overline{\mu},\overline{\nu}) = \mathcal{G}(\mu,\nu)$ follows upon observing that
	\begin{align*}
	\langle x^\alpha y^\beta Z^\gamma, \overline{\mu} \rangle
	&= 	\oneOverG \sum_{(A,B) \in \mathbb{G}}\langle x^\alpha y^\beta Z^\gamma, G_{A,B} \sharp \mu \rangle \nonumber \\[-0.35ex]
	&= 	\oneOverG \sum_{(A,B) \in \mathbb{G}}\langle (Ax)^\alpha (By)^\beta (BZA^\tr)^\gamma, \mu \rangle \nonumber \\[-0.35ex]
	&= 	\oneOverG \sum_{(A,B) \in \mathbb{G}}\langle x^\alpha y^\beta Z^\gamma, \mu \rangle \nonumber \\
	&= 	\langle x^\alpha y^\beta Z^\gamma, \mu \rangle
	\end{align*}
	for all $(\alpha,\beta,\gamma) \in \mathcal{S}$. The first equality is a direct consequence of the linearity of integration against measures; the second one follows from the definition of the pushforward measure; the third one exploits assumption~\ref{ass:invariance-4}.
	
	There remains to show that $\overline{\mu}$ and $\overline{\nu}$ satisfy~\eqref{e:invariance-mu}, \eqref{e:invariance-nu}, and~\eqref{e:meas-cond-div}--\eqref{e:meas-cond-bc}. For~\eqref{e:invariance-mu}, observe that
	\begin{align*}
	\langle \eta \circ G_{A,B}, \overline{\mu} \rangle
	&= \oneOverG \sum_{(A',B') \in \mathbb{G}} \langle \eta \circ G_{A,B}, G_{A',B'} \sharp \mu \rangle \nonumber \\[-0.35ex]
	&= \oneOverG \sum_{(A',B') \in \mathbb{G}} \langle \eta \circ G_{A,B} \circ G_{A',B'}, \mu \rangle \nonumber \\[-0.35ex]
	&= \oneOverG \sum_{(A'',B'') \in \mathbb{G}} \langle \eta \circ G_{A'',B''}, \mu \rangle \nonumber \\[-0.35ex]
	&= \oneOverG \sum_{(A'',B'') \in \mathbb{G}} \langle \eta, G_{A'',B''} \sharp \mu \rangle \nonumber \\
	&= \langle \eta, \overline{\mu} \rangle.
	\end{align*}
	The first and last equalities exploit the linearity of integration against measures; the second and fourth ones follow from the definition of the pushforward measure; the third one combines conditions~\ref{ass:group-2} and~\ref{ass:group-3} with the definition of $G_{A,B}$. Analogous steps show that $\overline{\nu}$ satisfies~\eqref{e:invariance-nu}.
	
	Next, we show that $\overline{\mu}$ and $\overline{\nu}$ satisfy~\eqref{e:meas-cond-div}.
	Fix an arbitrary function $\varphi \in \Phi$ and, for every $(A,B) \in \mathbb{G}$, define another function $\varphi_{A,B} \in \Phi$ via
	\begin{equation*}
	\varphi_{A,B}(x,y) := A^\tr \varphi(Ax, By).
	\end{equation*}
	Thanks to condition~\ref{ass:group-1}, it is an exercise in vector calculus to check that
	$\mathcal{D}\varphi_{A,B} = \mathcal{D}\varphi \circ G_{A,B}$.
	Consequently,
	\begin{align*}
	\langle \mathcal{D}\varphi, \overline{\mu} \rangle
	&= \oneOverG \sum_{(A,B) \in \mathbb{G}}  \langle \mathcal{D}\varphi, G_{A,B} \sharp \mu \rangle \nonumber \\[-0.35ex]
	&= \oneOverG \sum_{(A,B) \in \mathbb{G}}  \langle \mathcal{D}\varphi \circ G_{A,B}, \mu \rangle \nonumber \\[-0.35ex]
	&= \oneOverG \sum_{(A,B) \in \mathbb{G}}  \langle \mathcal{D}\varphi_{A,B}, \mu \rangle \nonumber \\[-0.35ex]
	&= \oneOverG \sum_{(A,B) \in \mathbb{G}}  \langle \varphi_{A,B} \cdot \hat{n}, \nu \rangle.
	\end{align*}
	The last equality is true because $(\mu,\nu) \in \mathcal{A}$ by assumption, so $\mu$ and $\nu$ satisfy~\eqref{e:meas-cond-div}. Now, assumption~\ref{ass:invariance-1} implies that the unit normal vector $\hat{n}(x)$ to $\partial\Omega$ satisfies $\hat{n}(Ax) = A \hat{n}(x)$, so $\varphi_{A,B}(x,y) \cdot \hat{n}(x) = (\varphi \cdot \hat{n}) \circ H_{A,B}$. Consequently,
	\begin{multline*}
	\langle \mathcal{D}\varphi, \overline{\mu} \rangle
	= \oneOverG \sum_{(A,B) \in \mathbb{G}}  \langle (\varphi \cdot \hat{n}) \circ H_{A,B}, \nu \rangle \\
	= \oneOverG \sum_{(A,B) \in \mathbb{G}}  \langle \varphi \cdot \hat{n} , H_{A,B} \sharp \nu \rangle 
	= \langle \varphi \cdot n , \overline{\nu} \rangle.
	\end{multline*}
	Since $\varphi \in \Phi$ was arbitrary, the pair $(\overline{\mu},\overline{\nu})$ satisfies~\eqref{e:meas-cond-div} as claimed. The remaining conditions, \eqref{e:meas-cond-marginal-mu}--\eqref{e:meas-cond-bc}, are verified with similar steps (omitted for brevity).
	\myqed
\end{pf}

\section{Reduction to an SDP}\label{s:sdp-projection}

The measure-theoretic inequality~\eqref{e:G-suff-condition} and its symmetry-exploiting analogue~\eqref{e:G-suff-condition-symm}, although simpler than the original functional inequality $\mathcal{F}(u) \geq 0$ from Problem~\ref{p:main-problem}, remain difficult to check in general. Often, however, they can be further strengthened into finite-dimensional sufficient conditions that can be tested computationally. One case in which this can be done systematically is that of problems with polynomial data, meaning that:
\begin{enumerate}[({P}1), leftmargin=*]
	\item\label{ass:poly-f} The function $f$ in~\eqref{e:polynomial-of-integrals} is polynomial;
	\item\label{ass:poly-ab} The functions $a$ and $b$ in~\eqref{e:pde-bcs} are polynomials;
	\item\label{ass:poly-domain} The domain $\Omega$ is a basic semialgebraic set, i.e., it is defined by polynomial equalities and inequalities.
\end{enumerate}
We shall also assume that:
\begin{enumerate}[({P}1), resume, leftmargin=*]
	\item\label{ass:poly-omega} 
		There exist polynomials $g$ and $h$ such that 
		\begin{subequations}
			\label{e:semialg-domain}
			\begin{align}
			\Omega &= \{x \in \mathbb{R}^n: g(x) \geq 0 \},\\
			\partial\Omega &= \{x \in \mathbb{R}^n: g(x) =  0 \}.
			\end{align}
		\end{subequations}
	\item\label{ass:poly-normal} The unit normal vector $\hat{n}(x) = \nabla g(x) \abs{\nabla g(x)}^{-1}$ to $\partial\Omega$ is polynomial.
\end{enumerate}
%

The last two assumptions are made purely for simplicity.
It is straightforward to extend our discussion to domains with semialgebraic definitions more complicated than~\eqref{e:semialg-domain}. When $\hat{n}(x)$ is not polynomial, instead, one can redefine the boundary measures in \S\ref{s:ob-measures} to be the pushforward of the weighted surface measure $\abs{\nabla g(x)}^{-1}\dSurf$. Then, all results in \S\ref{s:om-lifting} hold with $\hat{n}(x)$ replaced by $\nabla g(x)$ and with an extra factor of $\abs{\nabla g(x)}^{-1}$ in the surface integrals in~\eqref{e:bnd-measure-def} and~\eqref{e:meas-cond-marginal-nu}.

\subsection{A sufficient polynomial inequality}\label{s:poly-inequality}
Under assumptions~\ref{ass:poly-f}--\ref{ass:poly-normal}, the measure-theoretic inequalities~\eqref{e:G-suff-condition} and~\eqref{e:G-suff-condition-symm} can be strengthened into polynomial inequalities over sets defined by linear equations and linear matrix inequalities (LMIs). 
To do this, observe that $\mathcal{G}$ depends on the moments of $\mu$ with exponents $\mathcal{S} \subset \mathcal{E}_{d,p}$, which have degree $\abs{\alpha}\leq d$ and $\abs{\gamma}+\abs{p} \leq p$ (cf. \S\ref{s:problem-statement}). Recall also that the polynomials $a(x,y,Z)$ and $b(x,y)$ in~\eqref{e:pde-bcs} cannot have degree larger than $p$ in $y$ and $Z$, and there is no loss of generality in assuming that their degree in $x$ is $d$ or less. We may similarly assume that the polynomial $g$ in~\eqref{e:semialg-domain} has degree no larger than $d$.

Now, fix any integer $\omega \geq d/2$ (called the \textit{relaxation order}) and define the moment vectors
\begin{subequations}
	\begin{align*}
	\xi &:= \{ \langle x^\alpha y^\beta Z^\gamma, \mu \rangle: \; \abs{\alpha} \leq 2\omega,\; \abs{\beta}+\abs{\gamma} \leq p \},\\
	\theta &:= \{ \langle x^\alpha y^\beta, \nu \rangle: \; \abs{\alpha} \leq 2\omega + \deg(\hat{n}),\; \abs{\beta}\leq p  \}.
	\end{align*}
\end{subequations}
%
%
The functional $\mathcal{G}$ in~\eqref{e:G-def} is clearly a polynomial of $\xi$. If $2\omega=d$ we have exactly $\mathcal{G}(\mu,\nu) = f(\xi)$, where $f$ is the polynomial in~\eqref{e:polynomial-of-integrals}. With a slight abuse of notation we write the same when $2\omega > d$, which amounts to viewing $f$ as a higher-dimensional polynomial. Thus, to prove inequalities~\eqref{e:G-suff-condition} and~\eqref{e:G-suff-condition-symm} it suffices to show that $f$ is nonnegative for all moment vectors $\xi$. In turn, this is true if $f(\xi)\geq 0$ for all $\xi$ in a set defined by LMIs and linear equations necessarily satisfied by all moment vectors.

The LMIs involve the well-known \emph{moment matrices} and \emph{localizing matrices} \citep{Laurent2009,Lasserre2015}. Precisely, given any polynomial $h(x) = \sum_{\abs{\alpha}\leq d_h} h_\alpha x^\alpha$ of degree $d_h \leq d$, let $M_h(\xi)$ be a matrix with rows and columns indexed by the multi-indices $(\alpha,\beta,\gamma)$ with $\abs{\alpha}\leq \omega-\lceil \frac12 d_h \rceil$ and $\abs{\beta}+\abs{\gamma}\leq\lfloor \frac12 p \rfloor$, and whose entries are
\begin{equation*}
[M_h(\xi)]_{(\alpha,\beta,\gamma),(\alpha',\beta',\gamma')} = \sum_{\abs{\alpha''}\leq d_h} h_{\alpha''} \xi_{\alpha+\alpha'+\alpha'', \beta+\beta', \gamma+\gamma'}
\end{equation*}
Let $M_h(\theta)$ be defined in an analogous way, except that its rows and columns are indexed by multi-indices $(\alpha,\beta)$ with $\abs{\alpha}\leq \omega-\lceil \frac12 d_h \rceil$ and $\abs{\beta}\leq\lfloor \frac12 p \rfloor$. Then, under assumptions~\ref{ass:poly-f}--\ref{ass:poly-normal}, the moment vectors $\xi$ and $\theta$ are known to satisfy the block-diagonal LMI
\begin{equation}\label{e:lmi}
M(\xi,\theta) := \begin{pmatrix}
M_1(\xi)\\
& M_g(\xi)\\
&&M_1(\theta)
\end{pmatrix}
\succeq 0.
\end{equation}
%
%


The linear equality constraints on $\xi$ and $\theta$, instead, are obtained from~\eqref{e:meas-cond-div}--\eqref{e:meas-cond-bc} and (in the presence of symmetries) from~\eqref{e:invariance-mu} and~\eqref{e:invariance-nu}. This is done by restricting the test functions in these conditions to be monomials of suitably bounded degree. Let $\Phi_\omega \subset \Phi$ be a monomial basis for the space of $n$-valued polynomials $\varphi(x,y)$ of degree $2\omega$ in $x$ and degree $p$ in $y$:
%
\begin{equation*}
\Phi_{\omega} := \big\{ x^\alpha y^\beta \hat{e}_k: \;
\abs{\alpha} \leq 2\omega, \,
\abs{\beta}\leq p, \,
k=1,\ldots,n \big\},
\end{equation*}
where $\hat{e}_k$ is the $n$-dimensional unit vector in the $k$-th coordinate direction. Further, define the monomial sets 
%
\begin{subequations}
	\begin{align*}
	\Psi_{a,\omega} &:= \big\{ x^\alpha y^\beta Z^\gamma: \, 
	\abs{\beta} + \abs{\gamma} + \deg_{y,Z}(a) \leq p, \nonumber\\
	&\hspace{96pt}\abs{\alpha} + \deg_x(a) \leq 2\omega
	\big\},
	\\
	\Psi_{b,\omega} &:= \big\{ x^\alpha y^\beta: \, 
	\abs{\beta} + \deg_{y}(b) \leq p, \nonumber\\
	&\hspace{52pt}\abs{\alpha} + \deg_x(b) \leq 2\omega+\deg(\hat{n})
	\big\}.
	\end{align*}
\end{subequations}
Then, conditions~\eqref{e:meas-cond-div}--\eqref{e:meas-cond-bc} imply
\begin{subequations}
	\begin{align}
	\label{e:meas-cond-div-moment}
	\langle \mathcal{D}\varphi, \mu  \rangle &= \langle \varphi \cdot \hat{n} , \nu \rangle &&\forall \varphi \in \Phi_{\omega},\\
	\label{e:meas-cond-marginal-mu-moment}
	\langle x^\alpha, \mu \rangle &= \textstyle\int_{\Omega} x^\alpha \,\dVolume  &&\forall \alpha: \abs{\alpha} \leq 2\omega,\\
	\label{e:meas-cond-marginal-nu-moment}
	\langle x^\alpha, \nu \rangle &= \textstyle\int_{\partial\Omega} x^\alpha \,\dSurf &&\forall \alpha:  \abs{\alpha} \leq 2\omega+\deg(\hat{n}),\\
	\label{e:meas-cond-diff-moment}
	\langle \psi a, \mu \rangle &= 0 &&\forall \psi \in \Psi_{a,\omega},\\
	\label{e:meas-cond-bc-moment}
	\langle \rho b, \nu \rangle &= 0 &&\forall \rho \in \Psi_{b,\omega}.
	\end{align}
\end{subequations}
Thanks to the imposed degree bounds, these identities are linear equations for the moment vectors $\xi$ and~$\theta$.
For problems with symmetries, further linear equations can be obtained by setting $\eta = x^\alpha y^\beta Z^\gamma$ with $\abs{\alpha}\leq 2\omega$ and $\abs{\beta}+\abs{\gamma} \leq p$ in~\eqref{e:invariance-mu}, and by taking $\zeta = x^\alpha y^\beta$ with $\abs{\alpha}\leq 2\omega$ and $\abs{\beta}\leq p$ in and~\eqref{e:invariance-nu}. For notational convenience, we combine all equality constraints into the equation
\begin{equation}\label{e:equalities}
e(\xi,\theta) := A \xi + B \theta + c = 0.
\end{equation}
%
The matrices $A$ and $B$ and the vector $c$ can be constructed explicitly in applications.

The next result is immediate.
\begin{prop}\label{prop:poly-inequality}
	Inequalities~\eqref{e:G-suff-condition} and~\eqref{e:G-suff-condition-symm}, and hence the functional inequality in Problem~\ref{p:main-problem}, hold if
	\begin{equation}\label{e:final-poly-inequality}
	f(\xi) \geq 0 \quad \forall (\xi,\theta) \text{ s.t.~\eqref{e:lmi} and~\eqref{e:equalities}}.
	\end{equation}
\end{prop}

\begin{rem}
	The sufficient condition for the functional inequality in Problem~\ref{p:main-problem} provided by Proposition~\ref{prop:poly-inequality}  becomes less restrictive when the relaxation order $\omega$ (hence, the size of the vectors $\xi$ and $\theta$) is increased. However, in general we do not expect inequality~\eqref{e:final-poly-inequality} to become a necessary condition for Problem~\ref{p:main-problem} as $\omega \to \infty$. 
\end{rem}

\subsection{Implementation via semidefinite programming}\label{s:sos-relaxation}

Inequality~\eqref{e:final-poly-inequality} can be verified computationally via semidefinite programming using (a variation of) known techniques based on sum-of-squares (SOS) polynomials and SOS polynomial matrices \citep{Scherer2006}.

We begin by stating some preliminary results about positive semidefinite (PSD) matrices. For any integers $m$ and $n$, let $\mathbb{S}^n$ be the space of $n\times n$ symmetric matrices and let $\mathbb{S}^n_+$ be the PSD cone. Given $P \in \mathbb{S}^n$ and $Q \in \mathbb{S}^{mn}$ partitioned into $m^2$ blocks $Q_{ij}$ of size $n \times n$, define 
\begin{equation*}
\mathcal{T}(P,Q) := \begin{pmatrix}
\trace(PQ_{11}) & \cdots & \trace(PQ_{1n}) \\
\vdots & \ddots & \vdots \\
\trace(PQ_{n1}) & \cdots & \trace(PQ_{nn}) \\
\end{pmatrix}\in \mathbb{S}^m.
\end{equation*}
Note that $\mathcal{T}(P,Q)$ is simply the trace inner product of $P$ and $Q$ if $m=1$. When $m=n^k$ for some integer $k\geq1$, we define the matrices $\mathcal{T}_i(P,Q) \in \mathbb{S}^{n^{k-i}}$ for $i=1,\ldots,k$ via 
\begin{align*}
\mathcal{T}_1(P,Q) &= \mathcal{T}(P,Q), \\
\mathcal{T}_i(P,Q) &= \mathcal{T}[ P, \, \mathcal{T}_{i-1}(P,Q)], \quad i=2,\ldots,k.
\end{align*}
Observe that $\mathcal{T}_i(P,Q)$ depends linearly on $Q$ for every $i$, and that $\mathcal{T}_k(P,Q)$ is a scalar.

\begin{lem}\label{lemma:psd-matrices}
	If $P\in\mathbb{S}^n_+$ and $Q\in\mathbb{S}^{n^k}_+$, then $\mathcal{T}_i(P,Q)$ is PSD for each $i=1,\ldots,k$. In particular, $\mathcal{T}_k(P,Q) \geq 0$.
\end{lem}

\begin{pf}
	The case $i=1$ was proven by \citet[\S2.2]{Scherer2006}. The general case follows by induction. \myqed
\end{pf}

Lemma~\ref{lemma:psd-matrices} immediately implies the next result, which may be seen as an application of the so-called S-procedure.

\begin{prop}\label{prop:suff-condition-polynomial-s-procedure}
	Let $\rho_\omega$ be the number of rows in the LMI~\eqref{e:lmi} and $\sigma_\omega$ be the number of constraints in~\eqref{e:equalities}.
	The polynomial inequality~\eqref{e:final-poly-inequality}, hence the functional inequality in Problem~\ref{p:main-problem}, holds if there exist a polynomial vector $q(\xi,\theta)$ with $\sigma_\omega$ entries and PSD polynomial matrices 
	$\{Q_k(\xi,\theta)\}_{k=1,\ldots,K}$, each of size $(\rho_\omega)^k \times (\rho_\omega)^k$, such that
	%
	\begin{equation}\label{e:s-proc}
	f(\xi) \geq q(\xi,\theta) \cdot e(\xi,\theta) + \sum_{k=1}^K \mathcal{T}_k(M(\xi,\theta),Q_k(\xi,\theta)).
	\end{equation}
\end{prop}

It is well known (see, e.g., \citealp{Lasserre2015,Laurent2009,Scherer2006,Parrilo2013}) that if inequality \eqref{e:s-proc} and the positive semidefinitess constraints on the polynomial matrices $Q_k$ are strengthened into SOS constraints, then the search for these matrices and for the polynomial vector $q$ can be reformulated as an SDP once the number $K$ and the polynomial degrees are fixed. The same is true when $f$ depends affinely on parameters $\lambda$, as in Problem~\ref{p:main-opt-problem}, because in such cases~\eqref{e:s-proc} is jointly affine in $\lambda$, $q$, and the matrices $Q_k$. Consequently, one can use semidefinite programming to (attempt to) verify the functional inequality in Problem~\ref{p:main-problem} and to compute feasible $\lambda$ for the convex optimization problem~\eqref{e:opt-problem-general}.

\enlargethispage{1ex}
\begin{rem}\label{remark:K-values}
	For $K=1$, inequality~\eqref{e:s-proc} reduces to the sufficient conditions for polynomial inequalities considered by \cite{Scherer2006}. This is enough when the set of $(\xi,\theta)$ satisfying~\eqref{e:lmi} and~\eqref{e:equalities} is compact and satisfies the so-called Archimedean condition \citep[Theorem~1]{Scherer2006}, but not otherwise: as demonstrated in~\S\ref{ss:inhomog-ex}, in general it is necessary to use $K>1$. Beyond compactness, we are not aware of conditions ensuring that~\eqref{e:s-proc} is feasible if $K$ and the degrees of $q$ and each $Q_k$ are sufficiently large.
\end{rem}


\section{Examples}\label{s:examples}

We conclude by showing that Proposition~\ref{prop:suff-condition-polynomial-s-procedure} works well on two examples with known analytical solutions.

\subsection{A Jensen--Poincar\'e inequality}
Consider the problem of finding the minimum $\lambda$ such that
\begin{equation}\label{e:poincare-jensen}
\lambda \int_{-1}^1 u_x^2 \, \dVolume - \left( \int_{-1}^1 u \, \dVolume \right)^2 \geq 0
\end{equation}
for all functions $u:(-1,1)\to \mathbb{R}$ satisfying the boundary conditions $u(\pm1)=0$. This problem fits the framework of~\S\ref{s:problem-statement} with $\Omega = (-1,1)$, $n=m=1$, $p=2$, $a(x,u,\nabla u)\equiv0$ and $b(x,u)=u$.
Analytical solution using the calculus of variation shows that~\eqref{e:poincare-jensen} holds if $\lambda\geq\frac23$. Here, we verify this via the sufficient condition in Proposition~\ref{prop:suff-condition-polynomial-s-procedure}. We use $\omega = 1$ and take a semialgebraic definition of $\Omega$ as in~\eqref{e:semialg-domain} with $g(x)=1-x^2$.

Using the vector $\xi = \{\xi_{\alpha,\beta,\gamma}: \,\alpha\leq 2, \beta+\gamma\leq 2\}$, the polynomial inequality corresponding to~\eqref{e:poincare-jensen} is
\begin{equation}\label{e:poincare-jensen-poly-ineq}
f(\xi) = \lambda \xi_{0,0,2} - \xi_{0,1,0}^2 \geq 0.
\end{equation}

To determine the set over which this inequality is required to hold, recall from Remark~\ref{remark:dirichlet-bcs} that we can fix $\nu = \nu_0$ to be the boundary measure of the zero function, so
$\langle x^\alpha, \nu \rangle = 1^\alpha - (-1)^\alpha$ for all $\alpha$ and 
$\langle x^\alpha y^\beta, \nu \rangle = 0$ for all $\alpha$ and all $\beta > 0$.
%
With this simplification, identities~\eqref{e:meas-cond-div-moment}--\eqref{e:meas-cond-bc-moment} give equality constraints on the entries of $\xi$ alone. Specifically,~\eqref{e:meas-cond-div-moment} yields
\begin{align}\label{e:poincare-jensen-eq1}
	2-\xi_{0,0,0} &= 0, &
	\xi_{1,0,0} &= 0, &
	\tfrac23 - \xi_{2,0,0}&= 0.
\end{align}
Condition~\eqref{e:meas-cond-marginal-mu-moment}, instead, gives
\begin{subequations}
	\label{e:poincare-jensen-eq2}
	\begin{align}
		\xi_{0,0,1}  &= 0, &
		\xi_{0,1,1}  &= 0, \\
		2-\xi_{0,0,0} &= 0, &
		\xi_{0,1,0}+\xi_{1,0,1} &= 0,\\
		\xi_{0,2,0}+2\xi_{1,1,1} &= 0,&
		\xi_{1,0,0} &= 0,\\
		2\xi_{1,1,0}+\xi_{2,0,1} &= 0,&
		\xi_{1,2,0}+\xi_{2,1,1} &= 0.
	\end{align}	
\end{subequations}
Identities~\eqref{e:meas-cond-marginal-nu-moment}--\eqref{e:meas-cond-bc-moment} give no constraints because $\nu$ is fixed and $a=0$. However, since our example is invariant under the group $\mathbb{G}=\{(1,1), \, (-1,1)\}$ generating the variable transformation $(x,u,u_x) \mapsto (-x,u,-u_x)$ we can use~\eqref{e:invariance-mu} to obtain a final set of constraints:
\begin{subequations}
	\label{e:poincare-jensen-eq3}
	\begin{align}
	\xi_{0,0,1} &= 0,&
	\xi_{0,1,1} &= 0,&
	\xi_{1,0,0} &= 0,&
	\xi_{1,0,2} &= 0,\\
	\xi_{1,1,0} &= 0,&
	\xi_{1,2,0} &= 0,&
	\xi_{2,0,1} &= 0,&
	\xi_{2,1,1} &= 0.
	\end{align}	
\end{subequations}

The LMI~\eqref{e:lmi} for this example does not depend on $\theta$ and it is block-diagonal with two blocks,
\begin{equation}\label{e:poincare-jensen-lmi}
M(\xi,\theta) = \begin{pmatrix}
M_1(\xi) & 0 \\ 0 & M_g(\xi)
\end{pmatrix}\succeq 0.
\end{equation}
Using~\eqref{e:poincare-jensen-eq1}--\eqref{e:poincare-jensen-eq3} to eliminate 
variables that have a known value or can be expressed in terms of other ones, we find
	\begin{gather*}
	M_1(\xi)=
	{\small
		\begin{pmatrix}
	2 & \xi_{0,1,0} & -\xi_{0,1,0} & 0 & 0 & 0\\
	\xi_{0,1,0} & -2\xi_{1,1,1} & \xi_{1,1,1} & 0 & 0 & 0\\
	-\xi_{0,1,0} & \xi_{1,1,1} & \xi_{2,0,2} & 0 & 0 & 0\\
	0 & 0 & 0 & \xi_{0,0,2} & -\xi_{0,1,0} & \xi_{1,1,1}\\
	0 & 0 & 0 & -\xi_{0,1,0} & 2/3 & \xi_{2,1,0}\\
	0 & 0 & 0 & \xi_{1,1,1} & \xi_{2,1,0} & \xi_{2,2,0}		
	\end{pmatrix}
	}
	\intertext{and}
	M_g(\xi)=
	{\small
		\begin{pmatrix}
	4/3 &  \xi_{0,1,0}-\xi_{2,1,0} &   0  \\
	\xi_{0,1,0}-\xi_{2,1,0} &  -\xi_{1,1,1}-\xi_{2,2,0} &   0 \\
	0 &  0 &   \xi_{0,0,2}-\xi_{2,0,2}  
	\end{pmatrix}
	}.
	\end{gather*}

Next, we apply Proposition~\ref{prop:suff-condition-polynomial-s-procedure} to verify~\eqref{e:poincare-jensen-poly-ineq} for all vectors $\xi$ satisfying~\eqref{e:poincare-jensen-eq1}--\eqref{e:poincare-jensen-lmi}. Take $q(\xi,\theta)=0$, $K=1$, and $$Q_1(\xi,\theta) = \kappa^2 v v^\tr + \tfrac32 w(\xi) w(\xi)^\tr,$$ where $\kappa$ is an arbitrary constant and
\begin{align*}
v^\tr &= \left(0,\,0,\,0,\,1,\,0,\,0,\,0,\,0,\,0 \right)\\
w(\xi)^\tr &= \left(0,\,0,\,0,\,2/3,\,\xi_{0,1,0},\,0,\,0,\,0,\,0 \right).
\end{align*}
%
With these choices, inequality~\eqref{e:s-proc} reduces to
\begin{equation*}
\left(\lambda - \tfrac23 -\kappa^2 \right) \xi_{0,0,2} \geq 0
\end{equation*}
and it holds if $\lambda=\frac23+\kappa^2$. Since $\kappa$ is arbitrary, we conclude that the original functional inequality~\eqref{e:poincare-jensen} is satisfied for all $\lambda\geq \frac23$, as claimed at the start of the example.

\subsection{An inhomogeneous inequality}
\label{ss:inhomog-ex}
Let $\|v\|_2 = \int_{-1}^{1} v^2 \,\dVolume$ be the usual $L^2$ norm of a square-integrable function $v:(-1,1)\to\mathbb{R}$. Given any positive integer $d$, consider finding the minimum scalar $\lambda$ such that
\begin{equation}\label{e:inhomog-ex}
\lambda \|u_x\|_2^{2d} \|u\|_2^{2d} - 2 \|u\|_2^{2d} + 1 \geq 0
\end{equation}
for all $u:(-1,1)\to\mathbb{R}$ satisfying $u(\pm1)=0$. It follows from the Poincar\'e inequality $\|u_x\|_2^2 \geq (\pi/2)^2\|u\|_2^2$ that~\eqref{e:inhomog-ex} holds if $\lambda \geq (2/\pi)^{2d} =: \lambda_d^*$. 

Table~\ref{table:inhomog-ex} lists numerical upper bounds on $\lambda_d^*$ obtained for $d=1$, $2$, $3$ and $4$, computed using \texttt{MOSEK}~\citep{mosek} to solve the SDPs corresponding to the SOS strengthening of the sufficient condition~\eqref{e:s-proc} from Proposition~\ref{prop:suff-condition-polynomial-s-procedure}. This condition reads
\begin{equation}\label{e:inomog-ex-f}
f(\xi) := \lambda \xi_{0,0,2}^d \xi_{0,2,0}^d - 2 \xi_{0,2,0}^d + 1 \geq 0
\end{equation}
for all vectors $\xi = \{\xi_{\alpha,\beta,\gamma}:  \alpha\leq 2\omega,\,\beta+\gamma\leq 2 \}$ satisfying suitable equality constraints and an LMI $M(\xi)\succeq 0$, which can be explicitly constructed but are too long to be reported. As in our previous example, there is no dependence on $\theta$ by virtue of the Dirichlet boundary conditions on $u$. The tabulated bounds were obtained for increasing values of the relaxation order $\omega$ with $K=2$ and polynomial matrices $Q_1$ and $Q_2$ of degree $2d-2$, after enforcing all equality constraints on $\xi$ explicitly. Computations with $K=1$ give identical results if $d$ is even but are infeasible if $d$ is odd, showing that setting $K>1$ is necessary in general.

The bounds in Table~\ref{table:inhomog-ex} improve as the relaxation order $\omega$ is raised, and although they appear conservative they provably converge to the analytical optimum $\lambda^*_d$ as $\omega \to \infty$. Loosely speaking, this is because analysis by~\citet[Theorem~4.3]{Chernyavsky2021}, combined with a duality argument~\citep{Fantuzzi2019arxiv}, guarantees that our sufficient conditions capture with arbitrary accuracy the Poincar\'e inequality needed to verify~\eqref{e:inhomog-ex}. Precisely, for every $\kappa > (\pi/2)^2$ there exists a sufficiently large $\omega$ and a PSD matrix $P$ of the same size as the matrix $M(\xi)$ for that $\omega$ value such that, for all $\xi$ that satisfy the equality constraints,
\begin{equation*}
\xi_{0,0,2} = \kappa \xi_{0,2,0} + \trace[M(\xi) P].
\end{equation*}
Taking the $d$-th power of this identity, multiplying both sides by $\xi_{0,2,0}^d$, and rearranging gives
\begin{equation}\label{e:ex-inhomog-sum-expansion}
\xi_{0,0,2}^d \xi_{0,2,0}^d - \kappa^d \xi_{0,2,0}^{2d} = \sum_{j=1}^d \textstyle\binom{d}{j} \displaystyle\kappa^{d-j} \xi_{0,2,0}^{2d-j} s_j(\xi),
\end{equation}
where $s_j(\xi) = \left( \trace[M(\xi) P] \right)^j = (\xi_{0,0,2} - \kappa \xi_{0,2,0})^j$.
Now, it can be shown that $\xi_{0,2,0}$ is a diagonal entry of $M(\xi)$, so there exists a PSD diagonal matrix $D$ such that $\trace(DM(\xi)) = \xi_{0,2,0}$. Using this matrix, it is possible to rewrite the right-hand side of~\eqref{e:ex-inhomog-sum-expansion} as $\mathcal{T}_2[M(\xi),R(\xi)]$ with
\begin{multline*}
R(\xi) :=  \sum_{\substack{j=1\\j\text{ odd}}}^{d} \xi_{0,2,0}^{2d-j-1} s_{j-1}(\xi) \left( D\otimes P \right)\\[-3ex]
+ \sum_{\substack{j=2\\j\text{ even}}}^{d} \xi_{0,2,0}^{2d-j} s_{j-2}(\xi) \left( P\otimes P \right),
\end{multline*}
which can be easily shown to be an SOS matrix. Then, if $f(\xi)$ is as in~\eqref{e:inomog-ex-f} and we fix $K=2$, $q(\xi)=0$, $Q_1(\xi)=0$, and $Q_2(\xi)=\lambda R(\xi)$, the sufficient condition~\eqref{e:s-proc} from Proposition~\ref{prop:suff-condition-polynomial-s-procedure} for our example reduces to
\begin{align*}
0 
\leq f(\xi) - \mathcal{T}_2[M(\xi),Q_2(\xi)] 
= \lambda \kappa^d \xi_{0,2,0}^{2d}  - 2 \xi_{0,2,0}^d + 1.
\end{align*}
This inequality, hence the original functional inequality~\eqref{e:inhomog-ex}, holds if and only if $\lambda \geq \kappa^{-d}$. This construction works for every $\kappa > (\frac{\pi}{2})^2$ provided that $\omega$ is large enough, so the functional inequality~\eqref{e:inhomog-ex} can be verified with the sufficient condition in Proposition~\ref{prop:suff-condition-polynomial-s-procedure} for $\lambda$ as close to $(\frac{2}{\pi})^{2d}=\lambda_d^*$ as desired. 

\begin{table}
	\centering
	\label{table:inhomog-ex}
	\caption{Upper bounds on the minimum $\lambda$ for which inequality~\eqref{e:inhomog-ex} holds. The analytical minimum $\lambda_d^*$ is reported for comparison.}
	\begin{tabular}{ccccc}
		\toprule
		$\omega$ & $d=1$ & $d=2$ & $d=3$ & $d=4$\\
		\midrule
		1 & 4.000000 & 16.00000 & 64.00000 & 256.0000\\
		2 & 1.333333 & 1.777775 & 2.370369 & 3.160493\\
		3 & 0.854470 & 0.730117 & 0.623862 & 0.533069\\
		4 & 0.677588 & 0.459121 & 0.311073 & 0.210776\\
		5 & 0.590550 & 0.348746 & 0.205849 & 0.121443\\
		\midrule
		$\lambda_d^*$ & 0.405285 & 0.164256 & 0.066570 & 0.026980\\
		\bottomrule
	\end{tabular}
\end{table}


\begin{ack}
This work benefited from discussions with David Goluskin, Jason Bramburger, Alexandr Chernyavsky and Sergei Chernyshenko. 
\end{ack}

\bibliography{reflist-no-doi}

\end{document}